\documentclass[letterpaper, 10 pt, journal, twoside]{ieeetran}  

\IEEEoverridecommandlockouts                              

\usepackage{graphics}
\usepackage{graphicx}
 \newtheorem{thm}{Theorem}[section]
 \newtheorem{cor}[thm]{Corollary}
 {\rm}
 
  \newtheorem{rem}[thm]{Remark}
   \newtheorem{ex}[thm]{Example}
\usepackage{amsfonts,amssymb,amsopn,amscd}
 \usepackage{mathrsfs}
\def\x{\mathbf{x}}
\def\v{\mathbf{v}}

\def\o{\mathbf{\omega}}
\def\X{\mathbf{X}}

\def\M{\mathbf{M}}
\def\y{\mathbf{y}}
\def\z{\mathbf{z}}
\def\om{\mathbf{\Omega}}
\def\R{\mathbb{R}}
\def\N{\mathbb{N}}
\def\K{\mathbf{K}}
\def\P{\mathbf{P}}

\title{Representation of chance-constraints with strong
asymptotic guarantees}

\author{Jean B. Lasserre$^{1}$
\thanks{*This work was supported by the European Research Council (ERC) via an ERC-Advanced Grant for the \# 666981 project TAMING.
The author is very grateful to Tillmann Weisser for his help in the numerical experiments.}
\thanks{$^{1}$Jean B Lasserre is senior researcher at LAAS-CNRS and the Institute of Mathematics, University of Toulouse, France. 
        {\tt\small lasserre@laas.fr}}%
}

\begin{document}

\maketitle
\thispagestyle{empty}
\pagestyle{empty}

\begin{abstract}
Given $\epsilon \in (0,1)$, a probability measure $\mu$ on $\om\subset\R^p$ and a semi-algebraic set
$\K\subset\X\times\om$, we consider the feasible set 
$\X^*_\epsilon=\{\x\in\X:{\rm Prob}[(\x,\o)\in\K]\geq 1-\epsilon\}$
associated with  a chance-constraint. We provide a sequence
of outer approximations $\X^d_\epsilon=\{\x\in\X:h_d(\x)\geq0\}$, $d\in\N$, 
where $h_d$ is a polynomial of degree $d$ whose
vector of coefficients is an optimal solution of a semidefinite program.
The size of the latter increases with the degree $d$. We also obtain the strong 
and highly desirable asymptotic guarantee that $\lambda(\X^d_\epsilon\setminus\X^*_\epsilon)\to0$
as $d$ increases, where $\lambda$ is the Lebesgue measure on $\X$. Inner approximations with same guarantees
are also obtained.
\end{abstract}


\section{Introduction}

\IEEEPARstart{W} {e} consider the following general framework for decision under uncertainty : Let $\x\in\X\subset\R^n$ be a decision variable while $\o\in\R^p$ is a {\it disturbance} (or {\it noise}) 
parameter whose distribution $\mu$ 
(with support $\om\subset\R^p$) is known, i.e.,
its list of moments $\mu_\beta:=\int_\om  \o^\beta \,d\mu(\o)$, $\beta\in\N^p$,
is available in closed form or numerically.

Both $\x$ and $\o$  are linked by constraints of the form $(\x,\o)\in\K\subset\X\times\om$, where
\begin{equation}
\label{setk}
\K\,=\,\{\,(\x,\o):\:g_j(\x,\o)\,\geq\,0,\quad j=1,\ldots,m\},
\end{equation}
for some polynomials $(g_j)\subset\R[\x,\o]$, that is,
$\K$ is a basic semi-algebraic set.

Next, for each fixed $\x\in\X$, let $\K_\x\subset\om$ be the (possibly empty) set defined by:
\begin{equation}
\label{chance}
\K_\x\,:=\,\{\,\o\in \om:\: (\x,\o)\,\in\,\K\},\quad \x\in\X.
\end{equation}
Let $\epsilon\in (0,1)$ be fixed. 
The goal of this paper is to provide tight approximations of the set
\begin{eqnarray}
\label{chance1}
\X^*_\epsilon&:=&\{\,\x\in \X:\: \mu(\K_\x)\,\geq\,1-\epsilon\,\}\\
\nonumber
&=&\{\,\x\in \X:\: {\rm Prob}((\x,\o)\in\K)\,\geq\,1-\epsilon\,\}
\end{eqnarray}
in the form:
\begin{equation}
\label{chance2}
\X^d_\epsilon\,:=\,\{\,\x\in \X:\: h_d(\x)\,\geq\,0\,\},\qquad d\in\N,
\end{equation}
where $h_d$ is a polynomial of degree at most $d$. 

Such approximations are particularly useful for optimization 
and control problems with chance-constraints; for instance
problems of the form:
\begin{equation}
\label{chance3}
\min\,\{\,f(\x): \x\in C;\:{\rm Prob}((\x,\o)\in\K)\,\geq\,1-\epsilon\,\}.
\end{equation}
Indeed one then replaces problem (\ref{chance3}) with
\begin{equation}
\label{chance4}
\min\,\{\,f(\x): \x\in C;\: h_d(\x)\,\geq\,0\,\},
\end{equation}
where the uncertain parameter $\o$ has disappeared. So if $C$ is a basic semi-algebraic set
then (\ref{chance4}) is a standard polynomial optimization problem.
Of course the resulting decision problem (\ref{chance4}) may still be hard to solve because 
the sets $\X^*_\epsilon$ and $\X^d_\epsilon$ are not convex in general.
But this may be the price to pay for avoiding a too conservative formulation
of the problem.  The interested reader is referred to Henrion \cite{Henrion}, Pr\'ekopa \cite{Prekopa} and
Shapiro \cite{Shapiro} for a general overview of chance constraints in optimization and to
Calafiore and Dabbene \cite{control1}, Jasour et al. \cite{lagoa1} and Li et al. \cite{control2} in control (and the references therein).

However, in the formulation (\ref{chance4}) one has got rid of the disturbance parameter $\o$, 
and so one may 
apply the arsenal of Non Linear Programming algorithms to 
get a local minimizer of (\ref{chance4}). If $n$ is not too large
or if some sparsity is present in problem (\ref{chance4}) one may even run 
a hierarchy of semidefinite relaxations to approximate its global optimal value.
For the latter approach the interested reader is referred to \cite{book1} and
for a discussion on this approach to various control problems with chance constraints
we refer to the recent paper of Jasour et al. \cite{lagoa1} and the references therein.

In Jasour et al. \cite{lagoa1} the authors have considered some control problems
with chance constraints. They have provided an elegant formulation 
and a numerical scheme for solving the related problem of computing
\[\x^*\,=\,\arg\max\,\{\,\mu(\K_\x)\,:\: \x\in\X\,\}.\]
This problem is posed as an infinite-dimensional LP problem in an appropriate
space of measures, that is, a Generalized Moment Problem (GMP) as described in Lasserre \cite{book1}. 
Then to 
obtain $\x^*$ they solve a hierarchy of semidefinite relaxations, which is  the {\it moment-SOS} approach for solving the GMP. 
This GMP formulation has the particular and typical feature of including a constraint of {\it domination} $\phi\leq\psi$
 between two measures $\phi$ and $\psi$. Such domination constraints are particularly powerful and
 have been already used in a variety of different contexts. See for instance Henrion et al. \cite{sirev} 
 for approximating the Lebesgue volume of a compact semi-algebraic set, 
Lasserre \cite{gaussian} for computing Gaussian measures of semi-algebraic set, Lasserre \cite{lebesgue} for ``approximating" 
the Lebesgue decomposition of a measure with respect to another one. It has been used by Henrion and Korda \cite{korda1} 
for approximating regions of attraction, by  Korda et al. \cite{korda2} for approximating maximum controlled invariant sets,
and more recently in 
Jasour and Lagoa \cite{lagoa2}
for a unifying treatment of some control problems.

\subsection*{Contribution}

The approach that we propose for determining the set $\X^d_\epsilon$
defined in (\ref{chance2}) is very similar in spirit to that in \cite{sirev} and \cite{lagoa1} and can be viewed as an additional illustration of the versatility of the GMP
and the moment-SOS approach in control related problems.
Indeed we also define an 
infinite-dimensional LP problem $\P$ in an appropriate space of measures and a
sequence of semidefinite relaxations $(\P_d)_{d\in\N}$ of $\P$, whose associated monotone sequence of optimal values $(\rho_d)_{d\in\N}$ converges to the optimal value $\rho$ of $\P$. An optimal solution of the dual 
of $(\P_d)$ allows to obtain a polynomial $h_d$ of degree $2d$ whose super-level set
$\{\x:h_d(\x)\geq0\}$ is precisely the desired approximation $\X^d_\epsilon$ of $\X^*_\epsilon$ in (\ref{chance2}); in fact the sets
$(\X^d_\epsilon)_{d\in\N}$ provide a sequence of {\it outer} approximations of $\X^*_\epsilon$. 
We also provide the strong asymptotic guarantee that
\[\lim_{d\to\infty}\,\lambda(\X^d_\epsilon\setminus\X^*_\epsilon)\,=\,0,\]
where $\lambda$ is the Lebesgue measure on $\X$,
which to the best of our knowledge is the first result of this kind at this level of generality.
(The same methodology applied to chance-constraints of the form
${\rm Prob}((\x,\o)\in\K)<\epsilon$ would provide a sequence of {\it inner}  approximations
of the set $\{\x\in\X:{\rm Prob}((\x,\o)\in\K)<\epsilon\}$.)

Another contribution is to include a technique to accelerate the convergence
$\rho_d\to\rho$ which otherwise can be too slow. This technique is different from the one used in \cite{sirev}
for the related problem of computing the volume of a semi-algebraic set, and has the nice feature or preserving
the {\it monotonicity} of the convergence of $\rho_d\to\rho$. It can be applied whenever
$d\mu$ is the Lebesgue measure $d\o$ on $\om$ or $d\mu=q(\o)\,d\o$ or $d\mu=\exp(-q(\o))d\o$
for some homogeneous nonnegative polynomial $q$.

At last but not least, in principle we can also treat the case where
the support $\om$ of $\mu$ and the set $\K$ are not compact, which includes the important case where $\mu$ is the normal distribution. We briefly explain
what are the (technical) arguments which allow to extend the method to the 
non compact case.

Of course this methodology is computationally expensive and so far limited to relatively small size problems (but after all the problem is very hard). An interesting issue not discussed here is to investigate whether sparsity patterns can be exploited to handle problems with larger size.

\section{Notation, definitions and preliminary results}
\label{notation}
\subsection{Notation and definitions}
Let $\R[\x]$ be the ring of polynomials in the variables $\x=(x_1,\ldots,x_n)$ and let
$\R[\x]_d$ be the vector space of polynomials of degree at most $d$
(whose dimension is $s(d):={n+d\choose n}$).
For every $d\in\N$, let  $\N^n_d:=\{\alpha\in\N^n:\vert\alpha\vert \,(=\sum_i\alpha_i)\leq d\}$, 
and
let $\v_d(\x)=(\x^\alpha)$, $\alpha\in\N^n_d$, be the vector of monomials of the canonical basis 
$(\x^\alpha)$ of $\R[\x]_{d}$. 
A polynomial $f\in\R[\x]_d$ is written
\[\x\mapsto f(\x)\,=\,\sum_{\alpha\in\N^n}f_\alpha\,\x^\alpha.\]

Given a closed set $\mathcal{X}\subset\R^n$, denote by $M(\mathcal{X})$ the space of finite Borel measures on $\mathcal{X}$ and by $\mathcal{P}(\mathcal{X})$ the convex cone of polynomials that are nonnegative on $\mathcal{X}$.

\noindent
{\bf Moment matrix.} Given a sequence $\y=(y_\alpha)_{\alpha\in\N^n}$, let $L_\y:\R[\x]\to\R$ be the linear (Riesz) functional
\[f\:(=\sum_\alpha f_\alpha\, \x^\alpha)\:\mapsto L_\y(f)\,:=\,\sum_\alpha f_\alpha\,y_\alpha.\]
Given $\y$ and $d\in\N$, the {\it moment} matrix associated with $\y$, is the 
real symmetric  $s(d)\times s(d)$ matrix $\M_d(\y)$ with rows and columns indexed in $\N^n_d$ and with entries 
\[M_d(\y)(\alpha,\beta)\,:=\,L_\y(\x^{\alpha+\beta})\,=\,y_{\alpha+\beta},\quad \alpha,\beta\in\N^n_d.\]
{\bf Localizing matrix.} Given a sequence $\y=(y_\alpha)_{\alpha\in\N^n}$, and a polynomial $g\in\R[\x]$, the {\it localizing} moment matrix 
associated with $\y$ and $g$, is the 
real symmetric  $s(d)\times s(d)$ matrix $\M_d(g\,\y)$ with rows and columns indexed in $\N^n_d$ and with entries 
\begin{eqnarray*}
M_d(g\,\y)(\alpha,\beta)&:=&L_\y(g(\x)\,\x^{\alpha+\beta})\\
&=&\sum_\gamma g_\gamma\,y_{\alpha+\beta+\gamma},\quad \alpha,\beta\in\N^n_d.
\end{eqnarray*}

\subsection{The volume of a compact semi-algebraic set}

In this section we recall how to approximate as closely as desired the Lebesgue volume of a compact semi-algebraic set $\K\subset\R^n$. It will be the building block of the methodology to approximate the set $\X^*_\epsilon$ in (\ref{chance1}).

Let $\X\subset\R^n$ be a box and 
let $\lambda\in M(\X)$ be the Lebesgue measure on $\X$. Let $\K:=\{\x: g_j(\x)\geq0,\:j=1,\ldots,m\}$, assumed to be compact. For convenience and with no loss of generality we may and will assume that $g_1(\x)=M-\Vert\x\Vert^2$ for some $M>0$.
\begin{thm}[\cite{sirev}]
\label{th1}
Let $\K\subset\X$ and with nonempty interior. Then
\begin{equation}
\label{th1-eq1}
\lambda(\K)\,=\,\displaystyle\sup_{\phi\in M(\K)}\,\{\,\phi(\K):\:\phi\leq\lambda\},
\end{equation}
and $d\phi^*=1_\K(\x)d\lambda$ is the unique optimal solution.
\end{thm}
\vspace{0.2cm}
Problem (\ref{th1-eq1}) is an infinite-dimensional LP with dual
\begin{eqnarray}
\label{dual-f}
\rho&=&\displaystyle\inf_{p\in\mathcal{C}(\X)}\,\{\,\int_\X p\,d\lambda: \mbox{$p\geq0$ on $\X$; 
$p\geq1$  on $\K$}\}\\
\label{dual}
&=&\displaystyle\inf_{p\in\R[\x]}\,\{\,\int_\X p\,d\lambda: \mbox{$p\geq0$ on $\X$; 
$p\geq1$  on $\K$}\}
\end{eqnarray}
where $\mathcal{C}(\X)$ is the space of continuous functions on $\X$.
That (\ref{dual-f}) and (\ref{dual}) have same optimal value follows from compactness of $\X$
and Stone-Weierstrass Theorem. Next, as shown in \cite{sirev}, there is no duality gap, i.e., $\rho=\lambda(\K)$.  

\subsection{Semidefinite relaxations}

Let $d_j=\lceil {\rm deg}(g_j)/2\rceil$, $j=1,\ldots,m$. To approximate $\lambda(\K)$ one solves the hierarchy of semidefinite programs, indexed by $d\in\N$:
\begin{equation}
\label{primal}
\begin{array}{rl}
\rho_d\,=\,\displaystyle
\sup_{\y,\z}&\{\,L_\y(1): y_\alpha+z_\alpha=\lambda_\alpha,\quad \forall\alpha\in\N^n_{2d}\\
&\M_d(\y),\M_d(\z)\succeq0\\
&\M_{d-d_j}(g_j\,\y)\succeq0,\quad j=1,\ldots,m\}.
\end{array}
\end{equation}

\subsection*{Interpretation.}  Ideally, the variables $\y=(y_\alpha)$ (resp. $\z=(z_\alpha)$) of (\ref{primal}) should be viewed as ``moments" of the measure
$\phi$ (resp. the measure $\psi:=\lambda-\phi$) in (\ref{th1-eq1}) (and so $L_\y(1)=\phi(\K)$); the constraints $\M_d(\y),\M_{d-d_j}(g_j\,\y)\succeq0$ (resp. $\M_d(\z)\succeq0$) are 
precisely necessary conditions for the above statement to be true\footnote{As $\X=[-1,1]^n=\{\x: 1-x_j^2\geq0,\:j=1,\ldots,n\}$,
in principle one should also 
impose $\M_{d-1}((1-x_j^2)\,\z)\succeq0$, $j=1,\ldots,n$, for all $d$,
to ensure that $\psi$ is supported on $\X$. However as $\M_d(\z)\preceq \M_d(\lambda)$ for all $d$, in the limit as $d\to\infty$,
one has $\psi\leq\mu$ and so ${\rm support}(\psi)\subset \X$.}
(and which become sufficient as $d\to\infty$). 

The sequence $(\rho_d)_{d\in\N}$ is monotone non increasing and $\rho_d\to \lambda(\K)$ as $d\to\infty$. However the convergence
is rather slow and in \cite{sirev} the authors have proposed
to replace the criterion $L_\y(1)$ by $L_\y(h)$ where $h\in\R[\x]$ is a polynomial that is nonnegative on $\K$ and vanishes on the boundary of $\K$. If one 
denotes by $\y^d$ an optimal solution of (\ref{primal}) then
$\rho_d\to\int_\K hd\lambda$ and $y^d_0\to\lambda(\K)$ as $d\to\infty$. The convergence $y^d_0\to\lambda(\K)$ is much faster but is not monotone anymore, which can be annoying because we do not obtain a decreasing sequence of upper bounds on $\lambda(\K)$ as was the case with (\ref{primal}).
For more details the interested reader is referred to \cite{sirev}.

\subsection{Stokes can help}
\label{section-stokes}
This is why we propose another technique to accelerate the convergence
$\rho_d\to\lambda(\K)$ in (\ref{primal}) while maintaining its monotonicity.
So let $h\in\R[\x]$ be such that
$h(\x)=0$ for all $\x\in\partial\K$ (but $h$ is not required to be nonnegative on $\K$). Then by Stokes' theorem
(with vector field $X=e_i\in\R^n$, $e_{ij}=\delta_{i=j}$, $i,j=1,\ldots,n$), for each $\alpha\in\N^n$:
\[\int_\K\underbrace{\frac{\partial}{\partial x_i}(\x^\alpha\,h(\x))}_{\x\mapsto\theta^i_\alpha(\x)}\,d\lambda(\x)\,=\,0,\quad i=1,\ldots,n,\]
and so the optimal solution $\phi^*$ of Theorem \ref{th1} must satisfy
\[\int_\K\theta^i_\alpha(\x)\,d\phi^*(\x)\,=\,0,\quad \forall\alpha\in\N^n;\:i=1,\ldots,n.\]
Therefore in (\ref{primal}) we may impose the additional moment constraints
\begin{equation}
\label{stokes}
L_\y(\theta^i_\alpha)\,=\,0,\quad \forall\alpha\in\N^n_{2d-{\rm deg}(h)};\:i=1,\ldots,n.\end{equation}
To appreciate the impact of such additional constraints on the convergence
$\rho_d\to\lambda(\K)$, consider 
the simple example with $n=2$ and $\X=[-a,a]^2$, let $\K:=\{\x:\Vert\x\Vert^2\leq1\}$ so that
$\lambda(\K)=\pi$. For different values of $a$ and $d=3,4$, results are displayed in Table \ref{table2}.
~
\begin{table}[h]
\caption{The effect of Stokes constraints. $n=2$}
\begin{center}
\begin{tabular}{|c|c|c|c|c|}
\hline
 &$\rho_3$ &$\rho'_3$&$\rho_4$&$\rho'_4$\\
\hline
a =1.4  & 5.71&3.55 &5.38 &3.27\\
\hline
a=  1.3 & 5.38&3.41 &5.04 &3.21\\
\hline
a=   1.2&5.02 &3.31 &4.70 &3.18\\
\hline
a=   1.1&4.56 &3.25 &4.32 &3.16\\
\hline
a=   1.0&3.91 &3.20 &3.87 &3.15\\
\hline
\end{tabular}
\end{center}
\label{table2}
\end{table}
\begin{rem}
Theorem \ref{th1} is valid for any measure $\mu\in M(\X)$ and not only the Lebesgue measure $\lambda$. On the other hand,  additional Stokes constraints 
similar to (\ref{stokes}) (but now with vector field $X=\x$ and polynomial $\theta_\alpha$ below)
are valid provided that $d\mu=f\,d\lambda$ or $d\mu=\exp(f)\,d\lambda$ 
for some homogeneous polynomial $f\in\R[\x]$. Then with $d_f={\rm deg}(f)$,
\[\theta_\alpha(\x)\,=\,\left\{\begin{array}{ll}(n+d_f)\x^\alpha\,h+\langle \x,\nabla (\x^\alpha\,h)\rangle&\mbox{for $fd\lambda$}\\
\x^\alpha\,h(n+d_ff)+\langle\x,\nabla (\x^\alpha\,h)\rangle&\mbox{for $\exp(f)d\lambda$}.
\end{array}\right.\]
\end{rem}

\section{Main Result}
\label{main}
After the preliminary results  of \S \ref{notation}, we are now in position to state our main result. Let $\mu$ be the distribution of the noise parameter $\o\in\om$,
and let $\lambda$ be the Lebesgue measure on $\X$. The notation
$\lambda\otimes\mu$ denotes the product measure on $\X\times\om$, that is,
\[\lambda\otimes\mu(A\times B)\,=\,\lambda(A)\,\mu(B),\quad\forall
A\in\mathcal{B}(\X),\,B\in\mathcal{B}(\om).\]
With $\K\subset\X\times\om$ as in (\ref{setk}), and for every $\x\in\X$, let $\K_\x$ be as in (\ref{chance}) (possibly empty). Consider the infinite dimensional LP:
\begin{equation}
\label{chance-lp}
\rho\,=\,\displaystyle\sup_{\phi\in M(\K)}\,\{\,\phi(\K):\:\phi\leq\lambda\otimes\mu\}.
\end{equation}
\begin{thm}
\label{th2}
The unique optimal solution of (\ref{chance-lp}) is
\[d\phi^*((\x,\o))\,=\,1_\K((\x,\o))\,d\lambda\otimes\mu((\x,\o)),\]
and the optimal value $\rho$ of (\ref{chance-lp}) satisfies
\begin{eqnarray}
\nonumber
\rho&=&\int_{\X\times\om}1_\K((\x,\o))\,\lambda\otimes\mu(d(\x,\o))\\
\label{th2-eq2}
&=&\int_\X\mu(\K_\x)\,\lambda(d\x).
\end{eqnarray}
\end{thm}
\begin{IEEEproof}
That $\rho=\displaystyle\int_{\X\times\om}1_\K((\x,\o))\,\lambda\otimes\mu(d(\x,\o))$
follows from Theorem \ref{th1} (with $\lambda\otimes\mu$ instead of $\lambda$ in Theorem \ref{th1}). By Fubini-Tonelli's Theorem
(see Ash \cite{Ash}[Theorem 2.6.6, p. 103])
\[\int_{\X\times\om}1_\K((\x,\o))\,\lambda\otimes\mu(d(\x,\o))\,=\]
\[\int_\X\underbrace{\left(\int_\om1_\K((\x,\o))\,\mu(d\o)\right)}_{\mu(\K_\x)}\lambda(d\x)\,=\,
\int_\X\mu(\K_\x)\,\lambda(d\x)\]
\end{IEEEproof}

\subsection*{Semidefinite relaxations}
Let $d_j=\lceil{\rm deg}(g_j)/2\rceil$ for all $j$. As we did for (\ref{primal}) in \S \ref{notation},
let $\y=(y_{\alpha,\beta})$ and
$\z=(y_{\alpha,\beta})$, $(\alpha,\beta)\in\N^{n+p}$,
and relax (\ref{chance-lp}) to the following hierarchy
of semidefinite programs, indexed by $d\in\N$:
\begin{equation}
\label{chance-lp-relax}
\begin{array}{rl}
\rho_d\,=\,\displaystyle\sup_{\y,\z}&\{\,y_0:\\
\mbox{s.t.}&y_{\alpha,\beta}
+z_{\alpha,\beta}\,=\,\lambda_\alpha\cdot\mu_\beta,\quad (\alpha,\beta)\in\N^{n+p}_{2d}\\
&\M_d(\y),\,\M_d(\z)\succeq0\\
&\M_{d-d_j}(g_j\,\y)\,\succeq\,0,\quad j=1,\ldots,m\},
\end{array}
\end{equation}
and of course $\rho_{d}\geq\rho_{d+1}\geq\rho$ for all $d$.
The dual of (\ref{chance-lp-relax}) is the semidefinite program:
\begin{equation}
\label{chance-lp-dual}
\begin{array}{rl}
\rho^*_d\,=\,\displaystyle\inf_{p\in\R[\x,\o]_{2d}}&\{
\displaystyle\int_{\X\times\om}
p(\x,\o)\,\lambda\otimes\mu(d(\x,\o))\\
\mbox{s.t.}&p(\x,\o)\,\geq\,1,\quad\forall (\x,\o)\in \K\\
&\mbox{$p$ is SOS}\,\}.
\end{array}
\end{equation}
Again as $\K$ is compact, for technical reasons (but with no loss of generality) we may and will assume that in the definition (\ref{setk}) of $\K$, $g_1(\x)=M-\Vert\x\Vert^2$ for some $M>0$.
\begin{thm}
\label{th3}
Let $\K$ and $(\X\times\om)\setminus\K$ be with nonempty interior.
There is no duality gap between (\ref{chance-lp-relax}) and its dual (\ref{chance-lp-dual}), i.e., $\rho_d=\rho^*_d$ for all $d$. 
In addition (\ref{chance-lp-dual}) has an optimal solution $p^*_d\in\R[\x,\o]_{2d}$ such that
\[\rho_d=\rho^*_d=\int_{\X\times\om} p^*_d(\x,\o)\,\lambda\otimes \mu(d(\x,\o)).\]
Define $h^*_d\in\R[\x]_{2d}$ to be:
\[\x\mapsto h^*_d(\x)\,:=\,\int_\om p^*_d(\x,\o)\mu(d\o),\quad \x\in\R^n.\]
Then $h^*_d(\x)\geq\mu(\K_\x)$ for all $\x\in\X$ and
\[\rho_d=\int_\X h^*_d(\x)\,\lambda(d\x)\,\to\,\rho=\int_\X \mu(\K_\x)\,\lambda(d\x)\]
as $d\to\infty$.
\end{thm}
\begin{IEEEproof}
That $\rho_d=\rho^*_d$ is because Slater's condition holds for (\ref{chance-lp-relax}).
Indeed let $\y^*$ be the moments of $\phi^*$ in Theorem \ref{th2}
and $\z^*$ be the moments of $\lambda\otimes\mu-\phi^*$ (on
$(\X\times\om)\setminus\K$). Then as 
$\K$ has nonempty interior, $\M_d(\y^*)\succ0$ and $\M_d(g_j\,\y^*)\succ0$
for all $d$. Similarly as $(\X\times\om)\setminus\K$ has nonempty interior, 
$\M_d(\z^*))\succ0$. Moreover since the optimal value $\rho_d$ is finite for all $d$
this implies that (\ref{chance-lp-dual}) has an optimal solution $p^*_d\in\R[\x,\o]_{2d}$.
Therefore:
\begin{eqnarray*}
\rho_d&=&\displaystyle\int_{\X\times\om}
p^*_d(\x,\o)\,\lambda\otimes\mu(d(\x,\o))\\
&=&\displaystyle\int_\X(\underbrace{\int_\om p^*_d(\x,\o)\,\mu(d\o)}_{h^*_d(\x)\geq\mu(\K_\x)})\,\lambda(d\x)\,=\,
\int_\X h^*_d(\x)\,\lambda(d\x)
\end{eqnarray*}
where $h^*_d(\x)\geq\mu(\K_\x)$ follows from $p^*_d\geq1$ on $\K$. Finally
the convergence $\lim_{d\to\infty}\rho_d=\rho$  follows from Theorem \ref{th1}.
\end{IEEEproof}
Then as $h^*(\x)\geq\mu(\K_\x)$ on $\X$, the sets $\X^d_\epsilon=\{\x\in\X:h^*_d(\x)>1-\epsilon\}$, $d\in\N$, form a  
sequence of outer approximations of the set $\X^*_\epsilon$.
In fact more can be said.
\begin{cor}
\label{conv}
Let $h^*_d\in\R[\x,\o]_{2d}$ be as in Theorem \ref{th3}. 
Then the function $\x\mapsto \psi^*_d(\x):=h^*_d(\x)-\mu(\K_\x)$ is nonnegative on $\X$ and 
converges to $0$  in $L_1(\X,\lambda)$.  In particular 
$\psi^*_d\to 0$ in $\lambda$-measure\footnote{A sequence of functions $h,(h_n)_{n\in\N}$ on a measure space
$(\X,\mathcal{B},\lambda)$ converges in measure to $h$
if for all $\epsilon>0$, $\lambda(\{\x: \vert h_n(\x)-h(\x)\vert \geq\epsilon\})\to0$ as $n\to\infty$. The sequence converges almost-uniformly
to $h$ if to every $\epsilon>0$ there exists a set $B_\epsilon\subset\X$ such that $\lambda(B_\epsilon)<\epsilon$ and
$h_n\to h$ uniformly in $\X\setminus B_\epsilon$.}, and $\lambda$-almost uniformly
for some subsequence $(\psi^*_{d_k})_{k\in\N}$.
\end{cor}
\begin{IEEEproof}
As $\rho_d\to\rho$ as $d\to\infty$, 
\[\lim_{d\to\infty}\,\int_\X(\underbrace{h^*_d(\x)-\mu(\K_\x)}_{\geq0})\,\lambda(d\x)\,=\,0,\]
whence the convergence to $0$ in $L_1(\X,\lambda)$. Then convergence $\psi^*_d\to0$ in $\lambda$-measure,
and $\lambda$-almost sure
convergence for a subsequence follow from standard results from Real Analysis. See e.g. Ash \cite[Theorem 2.5.1]{Ash}.
\end{IEEEproof}

As we next see, the convergence $h^*_d(\x)\to\mu(\K_\x)$ in $\lambda$-measure established in Corollary \ref{conv} will be useful to obtain strong asymptotic guarantees.

\subsection{Strong asymptotic guarantees}

We here investigate asymptotic properties of 
the sequence of sets $(\X^d_\epsilon)_{d\in\N}$, as $d\to\infty$. 
\begin{cor}
\label{cor-final}
With $\X^*_\epsilon$ as in (\ref{chance1}), let $\X^d_\epsilon:=\{\x\in\X:h^*_d(\x)\geq 1-\epsilon\}$ where
$h^*_d$ is as in Theorem \ref{th3}, $d\in\N$. Then:
\begin{equation}
\label{cor-final-eq1}
\lim_{d\to\infty}\lambda(\X^d_\epsilon)\,=\,\lambda(\X^*_\epsilon).
\end{equation}
\end{cor}
\begin{IEEEproof}
Observe that
\[\X\setminus\X^*_\epsilon\,=\,\bigcup_{\ell=1}^\infty\{\x\in\X: \mu(\K_\x)<1-\epsilon-1/\ell\},\]
and therefore
\[\lambda(\X\setminus\X^*_\epsilon)\,=\,\lim_{\ell\to\infty}\lambda(\underbrace{\{\x\in\X: \mu(\K_\x)<1-\epsilon-1/\ell\}}_{R_\ell}).\]
Next, for each $\ell=1,\ldots$, write
\[\lambda(R_\ell)=\lambda(R_\ell\cap\{\x\in\X: h^*_d(\x)<1-\epsilon\})
+\lambda(R_\ell\cap \X^d_\epsilon).\]
By the convergence $h^*_d\to\mu(\K_\x)$ in $\lambda$-measure
as $d\to\infty$,
$\lim_{d\to\infty}\lambda(R_\ell\cap\X^d_\epsilon)=0$
and so
\begin{eqnarray*}
\lambda(R_\ell)&=&\lim_{d\to\infty}\lambda(R_\ell\cap\{\x\in\X: h^*_d(\x)\,<\,1-\epsilon\})\\
&\leq& \lim_{d\to\infty}\lambda(\{\x\in\X: h^*_d(\x)\,<\,1-\epsilon\})\\
&\leq&\lambda(\X\setminus\X^*_\epsilon).
\end{eqnarray*}
This implies
\[\lim_{d\to\infty}\lambda(\{\x\in\X: h^*_d(\x)\,<\,1-\epsilon\})=\lambda(\X\setminus\X^*_\epsilon),\]
which in turn yields the desired result (\ref{cor-final-eq1}).
\end{IEEEproof}

\subsection{Inner Approximations}

In the previous section we have provided a converging sequence 
$(\X^d_\epsilon)_{d\in\N}$ of outer approximations of $\X^*_\epsilon$.
Clearly, letting
$\K^c:=(\X\times\om)\setminus\K$,  the same methodology now applied to a chance constraint of the form
\[{\rm Prob}((\x,\o)\,\in\K^c)\,<\,\epsilon\]
would provide a converging sequence of {\it inner} approximations of the set
$\X^*_\epsilon:=\{\x\in\X: {\rm Prob}((\x,\o)\,\in\K)\,>\,1-\epsilon\}$.
To do so, (i) write $\K^c$ as a finite union $\bigcup_i\K^c_i$ of basic semi-algebraic sets $\K^c_i$
(whose $\mu\otimes\lambda$ measure of their overlaps is zero), (iii) apply the above methodology
to each $\K^c_i$, and then (ii) sum-up the results.

\subsection{Accelerating convergence}
\label{accelerate}
As we already have seen in Section \ref{section-stokes} for the semidefinite program
(\ref{primal}), as $d\to\infty$ the convergence 
$\rho_d\to\rho$ of the optimal value of (\ref{chance-lp-relax}) 
 can also be slow due to the Gibb's phenomenon\footnote{The Gibbs' phenomenon 
  appears at a jump discontinuity when one approximates a piecewise $C^1$ function with a continuous function, e.g. by its Fourier series.} that 
appears in the dual (\ref{chance-lp-dual}) when approximating the indicator function $\x\mapsto 1_\K(\x)$
by a polynomial. 

So assume that $\mu$ is the Lebesgue measure on $\om
$ where for instance $\om=[-1,1]^p$,
scaled to be a probability measure (but the same idea works
if $d\mu=h(\o)d\o$, or if $d\mu=\exp(h(\o))d\o$ for s ome homogeneous polynomial $h$).
Then again we propose to include additional constraints on the moments
$\y$ and $\z$ in (\ref{chance-lp-relax}) coming from additional properties of the optimal solution 
$\phi^*$ and $\psi^*=\lambda\otimes\mu-\phi^*$ of (\ref{chance-lp}). Again these additional properties are coming from Stokes' formula but now for integrals on $\K_\x$ (resp. $\om\setminus\K_\x$), then integrated over $\X$.

Let  $f_1,f_2\in\R[\x,\o]$ be the polynomials $(\x,\o)\mapsto f_1(\x,\o):=\prod_jg_j(\x,\o)$
and $(\x,\o)\mapsto f_2(\x,\o):=f_1(\x,\o)\prod_j (1-\o_j^2)$ of respective degree $d_1$, $d_2$. 
For each fixed $\x\in\X$, the polynomial
$\o\mapsto f_1(\x,\o)$ (resp. $\o\mapsto f_2(\x,\o)$)
vanishes on the boundary $\partial\K_\x$ of $\K_\x$ (resp. $\partial\K_\x^c$ of $\om\setminus\partial\K_\x$).
Therefore for each $\beta\in\N^p$,
Stokes' Theorem (applied with vector fields $e_j\in\R^p$ (where $e_{jk}=(\delta_{j=k})$), $k,j=1,\ldots,p$), states:
\[\displaystyle\int_{\K_\x} \frac{\partial}{\partial \o_j} (\o^\beta\,f_1(\x,\o))\,d\o=0,\quad \beta\in\N^p,\,j=1,\ldots,p,\]
\[\displaystyle\int_{\K_\x^c} \frac{\partial}{\partial\o_j} (\o^\beta\,f_2(\x,\o))\,d\o=0,\quad \beta\in\N^p,\,j=1,\ldots,p.\]
So let $\theta^k_{j\beta}\in\R[\x,\o]$ of degree $d_k+\vert\beta\vert-1$, $k=1,2$, be:
\[(\x,\o)\mapsto \theta^k_{j\beta}(\x,\o):=
\partial (\o^\beta\,f_k(\x,\o))/\partial\o_j,\]
for all $\beta\in\N^p$, $j=1,\ldots,p$.
Then for each $(\alpha,\beta)\in \N^{n+p}$:
\[\int_\X\int_{\K_\x}\x^\alpha\,\theta^1_{j\beta}(\x,\o)\,d\mu(\o)\,d\lambda(\x)\,=\,0,\]
\[\int_\X\int_{\K^c_\x}\x^\alpha\,\theta^2_{j\beta}(\x,\o)\,d\mu(\o)\,d\lambda(\x)\,=\,0.\]
Equivalently, in view of what are $\phi^*,\psi^*$ in Theorem \ref{th2},
\begin{eqnarray*}
\int_\K\x^\alpha\,\theta^1_{j\beta}(\x,\o)\,d\phi^*((\x,\o))&=&0,\quad(\alpha,\beta)\in\N^{n+p}\\
\int_{\K^c}\x^\alpha\,\theta^2_{j\beta}(\x,\o)\,d\psi^*((\x,\o))&=&0,\quad(\alpha,\beta)\in\N^{n+p}.\end{eqnarray*}
Therefore in (\ref{chance-lp-relax}) we may include the additional moments constraints
$L_\y(\x^\alpha\,\theta^1_{j\beta}(\x,\o))=0$, and $L_\z(\x^\alpha\,\theta^2_{j\beta}(\x,\o))=0$, 
for all $(\alpha,\beta)\in\N^{n+p}$ such that
$\vert\alpha+\beta\vert\leq 2d+1-d_1$ and 
$\vert\alpha+\beta\vert\leq 2d+1-d_2$ respectively.

\subsection{The non-compact case}

In some applications the noise $\o$ is assumed to follow a normal distribution 
$\mu$ on $\om=\R^p$. Therefore $\om$ is not compact anymore and the machinery used in \cite{sirev} cannot be applied directly. However 
the normal distribution satisfies the important Carleman's
property. That is, let $L_\y$ be the Riesz functional associated with $\mu$, 
i.e., $L_\y(f)=\int fd\mu$ for all $f\in\R[\o]$. Then
\begin{equation}
\label{carleman}
\sum_{k=1}^\infty L_\y(\o_i^{2k})^{-1/2k}\,=\,+\infty,\quad i=1,\ldots,p.
\end{equation}
In particular $\mu$ is {\it moment determinate}, that is, $\mu$
is completely determined by its moments.
These two properties have been used extensively in e.g. Lasserre \cite{lebesgue} and also in \cite{gaussian}, precisely to show that with $\K\subset\om$ not necessarily compact, one may still approximate its Gaussian measure $\mu(\K)$ as closely as desired. Again one solves (\ref{th1-eq1}) 
via the same hierarchy of semidefinite 
relaxations (\ref{primal}) (but now with $\mu$ instead of $\lambda$). 
For more details the interested reader is referred to Lasserre \cite{lebesgue,gaussian}. 

In view of the above (technical) remarks, one may then extend the machinery 
described in \S \ref{main} to the case where $\om=\R^p$, $\mu$ is the Gaussian measure, and $\K_\x$ ($\x\in\X$) is not necessarily compact.
A version of Stokes' Theorem for non compact sets 
is even described in \cite{gaussian} to accelerate the convergence of the semidefinite relaxations (\ref{primal}) (with $\mu$ instead of $\lambda$). 
It can be used to accelerate the convergence of the semidefinite relaxations (\ref{chance-lp-relax}), exactly as we do in \S \ref{accelerate} for the compact case.

\subsection{Numerical examples}
For illustration purposes we have considered simple small dimensional examples for which the function
$\x\mapsto\mu(\K_\x)$ has a closed form expression, so that we can compare the set $\X^*_\epsilon$
with its approximations $\X^d_\epsilon$, $d\in\N$, obtained in Corollary \ref{cor-final}
(with and without using Stokes constraints).

\begin{ex}
\label{ex3}
$\X=[-1,1]$, $\om=[0,1]$ and $\K=\{(x,\omega): 1-x^2/0.81-\omega^2/1.44\geq0\}$.
$\lambda$ and $\mu$ are the Lebesgue measure. In this case 
$\X^*_\epsilon=[a_\epsilon,b_\epsilon], \X^d_\epsilon=[a^d_\epsilon,b^d_\epsilon]\subset [-1,1]$. 
In Table \ref{table3} we display the relative error $(b^d_\epsilon-b_\epsilon+a_\epsilon-a^d_\epsilon)/
(b^d_\epsilon-a^d_\epsilon)$ for different values of $\epsilon$ and $d$, with and without Stokes constraints.
The results indicate that adding Stokes constraints help a lot. With relatively few moments $2d\leq 16$
one obtains good approximations.
\begin{table}[h]
\caption{Example \ref{ex3}; n=1: with and without Stokes \label{table3}}
\begin{center}
\begin{tabular}{|c|c|c|c|c|}
\hline
& d=4  &d=4 (Stokes) &d=8&d=8 (Stokes)\\
\hline
$\epsilon=0.75$ &13\% &3.8\% & 6.3\%& 0.7\%\\
\hline
$\epsilon=0.50$ &16.6\% & 2.1\% & 9.8\%& 0.2\%\\
\hline
$\epsilon=0.25$ &26.6\% &6.5\% & 18.5\%& 1.0\%\\
\hline
$\epsilon=0.00$ &44.7\% &22.7\%& 31.2\%&8.2\%\\
\hline
\end{tabular}
\end{center}
\end{table}
\end{ex}
\begin{ex}
\label{ex4}
$\X=[-1,1]$, $\om=[0,1]$ and $\K=\{(x,\omega): 1-x^2-\frac{\omega^2}{2}\geq0;\: \frac{x^2}{2}+\omega^2-\frac{1}{4}\geq0\}$.
$\lambda$ and $\mu$ are the Lebesgue measure on $\X$ and $\om$ respectively. In this case,
when $\epsilon<0.4$, the set  $\X^*_\epsilon$ is the union of two disconnected intervals, hence more difficult to approximate.
As in Example \ref{ex3}, in Table \ref{table4} we display the relative error 
for different values of $\epsilon$ and $d$, with and without Stokes constraints.
Again, the results indicate that adding the Stokes constraints help a lot. With relatively few moments $2d\leq 20$
one obtains good approximations. For instance,  $\X^*_{0.1}=[-0.7714,-0.3082]\cup [0.3082,0.7714]$
while one obtains  $\X^*_{0.1}\subset\X^{10}_{0.1}=[-0.7985,-0.26]\cup [0.26,0.7985]$ with Stokes and
the larger set $\X^{10}_{0.1}=[-0.8673,-0.1881]\cup [0.1881,0.8673]$ without Stokes constraints.
\begin{table}[h]
\caption{Example \ref{ex4}, n=1; with and without Stokes \label{table4}}
\begin{center}
\begin{tabular}{|c|c|c|c|c|}
\hline
& d=7  &d=7 (Stokes) &d=10&d=10 (Stokes)\\
\hline
$\epsilon=0.7$ &2.3\% &0.4\% & 2.3\%& 0.2\%\\
\hline
$\epsilon=0.4$ &30\% & 10\% & 26\%& 5.5\%\\
\hline
$\epsilon=0.1$ &52\% &24\% & 46\%& 16\%\\
\hline
\end{tabular}
\end{center}
\end{table}
\end{ex}
\begin{ex}
\label{ex5}
$\X=[-1,1]^2$, $\om=[0,1]$ and $\K=\{(x,\omega): 1-x_1^2-x_2^2-\omega^2\geq0\}$.
$\lambda$ and $\mu$ are the Lebesgue measure on $\X$ and $\om$ respectively. 
For this two-dimensional example (in $\x$) 
we have plotted the boundary of the set $\X^*_\epsilon$ (inner approximate circle, solid line in black). The 
curve in the middle (red dashed line) (resp. outer, blued dashed line) is the boundary of the approximation $\X^d_\epsilon$ computed with Stokes constraints
(resp. without Stokes constraints).
For $\epsilon=0.01$ and $d=10$
the results are plotted  in Fig. \ref{fig1} and in Fig. \ref{fig2}
for $\epsilon=0.05$ and $d=10$.

\begin{figure}
\begin{center}
\caption{\label{fig1}}
\includegraphics[width=.5\textwidth]{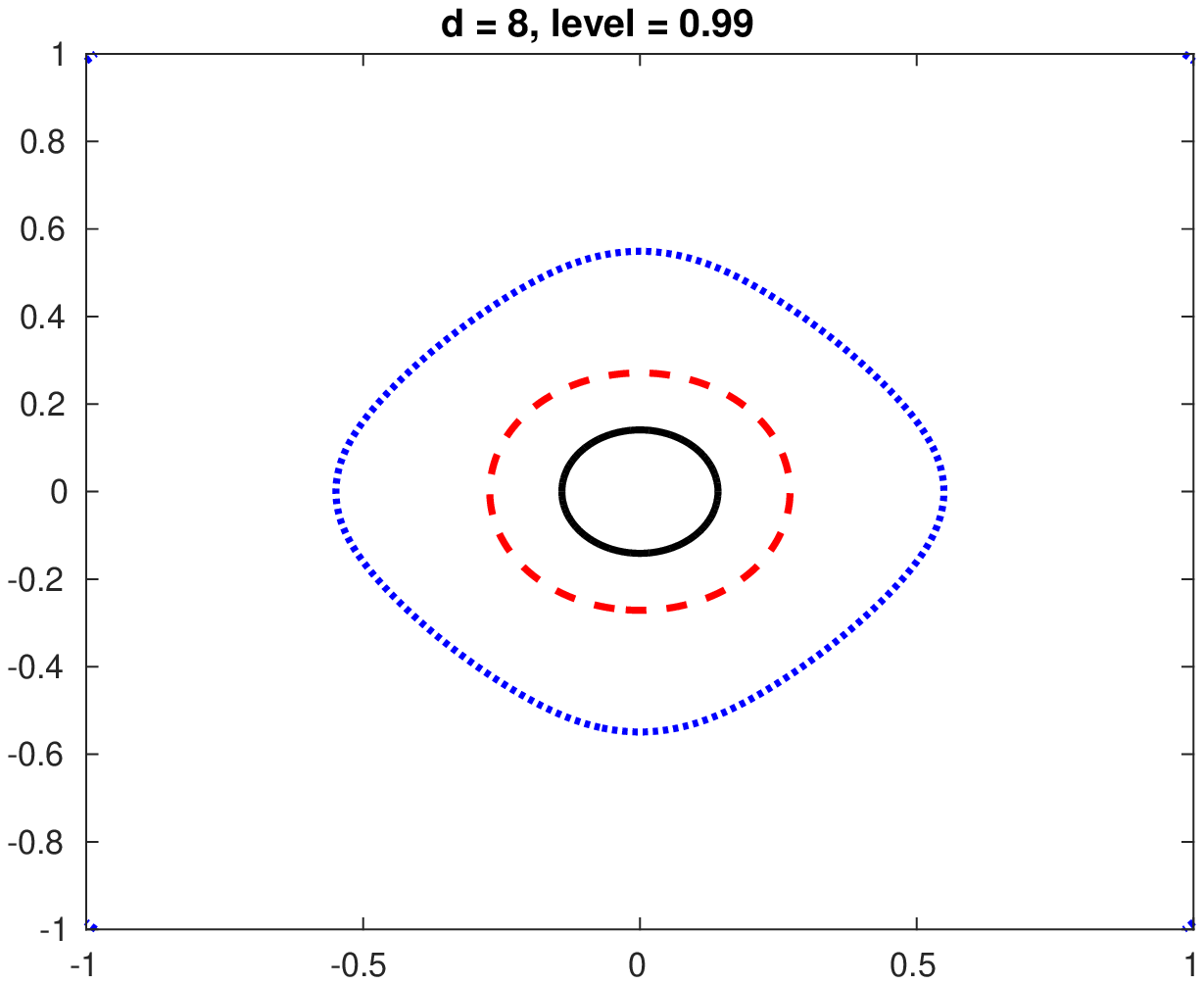}
\end{center}
\end{figure}

\begin{figure}
\begin{center}
\caption{ \label{fig2}}
\includegraphics[width=.5\textwidth]{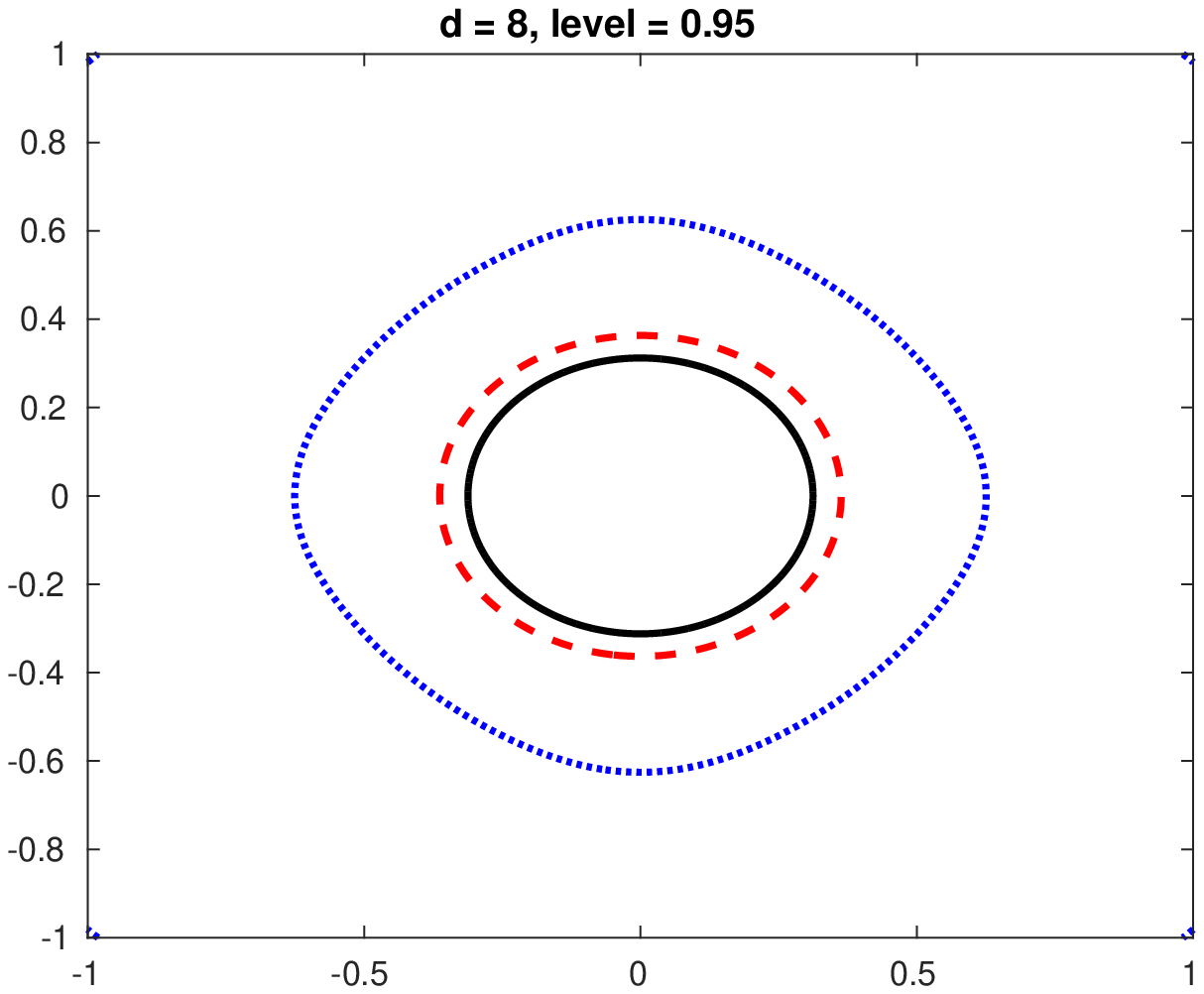}
\end{center}
\end{figure}
\end{ex}
\section{CONCLUSION}

We have presented a systematic numerical scheme to provide an
sequence of outer and inner approximations $(\X^d_\epsilon)$ of the feasible set $\X^*_\epsilon$  associated with 
chance constraints of the form ${\rm Prob}((\x,\o)\in\K)>1-\epsilon$.
Each outer and inner approximation $\X^d_\epsilon$ is the $0$-super level set of some polynomial whose coefficients are computed via solving a certain semidefinite program. As $d$ increases $\lambda(\X^d_\epsilon\setminus\X^*_\epsilon)\to 0$, a nice and highly desirable asymptotic property.  Of course this methodology
is computationally expensive and in its present form limited to problems of small size. But we hope it can pave the way to define efficient heuristics. 
Also checking whether
this methodology can accommodate potential sparsity patterns present in larger size problems, is a topic of further investigation.

\end{document}